# Shifted Monic Ultraspherical Approximation for solving some of Fractional Orders Differential Equations

**M. Abdelhakem[1,2]**
E-mail:mabdelhakem@yahoo.com

**Doha M.R.[1]**
E-mail:dohamahmoud147@gmail.com

**A. F. Saadallah[1]**
E-mail: adelsaadallah@hotmail.com

**M. El-Kady[1,2]**
E-mail: mamdouh_elkady@cic-cairo.com

[1]Mathematics Department, Faculty of Science, Helwan University, Cairo, Egypt.
[2]Canadian International College - CIC

**ABSTRACT**
The purpose of this paper is to show and explain a new formula that indicates with finality the derivatives of Shifted Monic Ultraspherical polynomials (SMUPs) of any degree and for any fractional-order using the shifted Monic Ultraspherical polynomials themselves. We also create a direct method solution for the linear or nonlinear multi-order fractional differential equations (FDEs) with constant coefficients involving a spectral Galerkin method. The spatial approximation with its fractional order derivatives (described in the Caputo sense) are built using shifted Monic Ultraspherical polynomials $\check{\mathbb{C}}_{L,j}^{\lambda}(\xi)$ with $\xi \in [0, L], L > 0$ and n is the polynomial degree.

**Introduction**

One of the old fields of mathematics is fractional calculus which dates back to the time of Leibniz (Abdeljawad, 2015) and from then many studies were done in this field (Akkurt et al., 2017; Garrappa and Popolizio, 2011). Fractional differential equations (FDEs) have attracted the interest of researchers in many areas such as Physics, Chemistry, Engineering and Social Sciences (Ortigueira *et al.,* 2015; Kilbas *et al.,* 2006). The analytic results on the existence and uniqueness of solutions to the FDEs have been investigated by many authors (El-Gamel *et al.,* 2003; Ortigueira *et al.,* 2015; Yang, 2013). Generally, most of the FDEs do not have analytic solutions, so one has to resort to approximation and numerical methods.

One class of FDEs is multi-order fractional differential equations. Since the last decade, extensive research has been conducted on the development of numerical methods for FDEs. In Li, (2010), Doha *et al.* (2011) proposed an efficient spectral tau and collocation method based on the Chebyshev polynomials for solving this equation. In dibi and Assari, (2010) the FDE is converted into a system of FDEs and the shifted Chebyshev operational matrix method is used to solve the resultant system. Some other works on this problem are: Haar wavelet method (Odibat and Momani, 2008) Legender wavelet method (Balaji, 2015) and second kind Chebyshev wavelet method dibi and Assari, (2010); Shiralashetti and Deshi, 2017).

In (Abdelhakem *et al.,* 2019; El-Kady and Moussa, 2011, 2013) the authors stated the formula of monic Chebyshev polynomials and used it to solve optimal control problems integral, integro-differential, and higher order differential equations. Also, in (Abdelhakem *et al.,* 2013) the authors used monic Ultraspherical to solve differential equations. Numerical examples are solved to show good ability and accuracy of the monic Chebyshev polynomials.

A shifted Monic Ultraspherical collation strategy is to be exhibited with SMU with equal step points to use as collation nodes for approaching multi-order fractional problems. Multiple numerical

**Corresponding Author:** Doha M. Radi, Faculty of Science, Helwan University, Cairo, Egypt.
E-mail:dohamahmoud147@gmail.com

examples are used to show the validity and suitability of the proposed techniques and to compare with the previous results. This method supported by examples of HBVPs in wide application. The mentioned examples showed that the proposed method is efficient and accurate.

The paper is coordinated as follows: In Section 2, we submit some needed definitions and give some significant properties of Ultraspherical polynomials. In Section 3, we state and prove the main result of the paper which gives undeniably a formula that asserts the fractional-order derivatives of the shifted Monic Ultraspherical polynomials in terms of the shifted Monic Ultraspherical polynomials themselves. In Section 4, we build and develop algorithms for solving FDEs by employing Galerkin spectral methods. In Section 5, there is a conclusion

## 1. Preliminaries and notation:

The fractional derivative in the Caputo sense
The Caputo's fractional differential is a linear operation

$$D^q(\gamma u(x) + \eta w(x)) = \gamma D^q u(x) + \eta D^q w(x), \quad (1)$$

Where $\gamma$ and $\eta$ are constant. For Caputo derivative we have

$$D^q C = 0, \quad (C \text{ is constant}). \quad (2)$$

$$D^q x^\beta = \begin{cases} 0, & \text{for } \beta \in N_0 \text{ and } \beta < \lceil q \rceil, \\ \frac{\Gamma(\beta+1)}{\Gamma(\beta+1-q)} x^{\beta-q}, & \text{for } \beta \in N_0 \text{ and } \beta \geq \lceil q \rceil \text{ or } \beta \notin N \text{ and } \beta > \lfloor q \rfloor. \end{cases} \quad (3)$$

We use the ceiling function $\lceil q \rceil$ to define as the smallest integer greater than or equal to $q$, and the floor function $\lfloor q \rfloor$ to define as the largest integer less than or equal to $q$. Also $N = \{1, 2, \dots\}$ and $N_0 = \{0, 1, 2, \dots\}$.

## 2. Properties of shifted Monic Ultraspherical polynomials:

In (Abdelhakem *et al.*, 2013) the authors use the analytic form of the Monic Ultraspherical polynomials. Also, Monic Ultraspherical polynomials are defined on the interval $[-1,1]$ and can be resolved with the aid of the following recurrence formula:

The essential recurrent formulae for Monic Ultraspherical polynomials are defined as: (4)

$$(j+1)\frac{2^j \Gamma(j+\lambda)}{\Gamma(\lambda)\Gamma(j+1)}\mathbb{C}_{j+1}^{(\lambda)}(t) = 2^{j+1}\frac{\Gamma(j+\lambda+1)}{\Gamma(\lambda)\Gamma(j+1)} t\mathbb{C}_j^{(\lambda)}(t) - (2\lambda+j-1)\frac{2^j \Gamma(j+\lambda)}{\Gamma(\lambda)\Gamma(j+1)}\mathbb{C}_{j-1}^{(\lambda)}(t) \quad (4)$$

where $\mathbb{C}_0^{(\lambda)}(t) = 1$ and $\mathbb{C}_1^{(\lambda)}(t) = t$. In order to use these polynomials on the interval $x \in [0, L]$ we defined the so-called Monic Shifted Ultraspherical polynomials by proposing the change of variable $t = \frac{2x}{L} - 1$. Let the Monic Shifted Ultraspherical polynomials $\mathbb{C}_j^{(\lambda)}\left(\frac{2x}{L} - 1\right)$ be denoted by $\check{\mathbb{C}}_j^{(\lambda)}(x)$. Then $\check{\mathbb{C}}_j^{(\lambda)}(x)$ can be created by using the following recurrence relation:

$$\check{\mathbb{C}}_{j+1}^{(\lambda)}(x) = \left(x - \frac{1}{2}\right)\check{\mathbb{C}}_j^{(\lambda)}(x) - \frac{(2\lambda+j-1)}{2^{j+2}(\lambda+j)(\lambda+j-1)}\check{\mathbb{C}}_{j-1}^{(\lambda)}(x), \quad j = 1,2,\cdots (5)$$

The analytic form of the MSUPs $\check{\mathbb{C}}_j^{(\lambda)}(x)$ of degree $j$ is given by

$$\check{\mathbb{C}}_j^{(\lambda)}(x) = \mathcal{H}_j \sum_{r=0}^{[j/2]} \sum_{k=0}^{j-2r} (-1)^{j-r-k} \frac{\Gamma(j-r+\lambda)}{\Gamma(\lambda)\Gamma(r+1)} \frac{2^{j-2r+k}}{\Gamma(k+1)\Gamma(j-2r-k+1)} x^k \quad (6)$$

where $[j/2]$ refers to the integer part of the fraction, $\mathcal{H}_j = 2^{-2j} \frac{\Gamma(\lambda) j!}{\Gamma(j+\lambda)}$, where $\check{\mathbb{C}}_0^{(\lambda)}(x) = 1$ and $\check{\mathbb{C}}_1^{(\lambda)}(x) = x - \frac{1}{4\lambda}$.

The orthogonally is

$$\int_0^1 (x-x^2)^{\lambda-\frac{1}{2}} \check{\mathbb{C}}_j^{(\lambda)}(x)\check{\mathbb{C}}_k^{(\lambda)}(x) dx = \mathcal{H}_j \mathcal{H}_k \delta_{jk} \psi_j \quad (7)$$

where $\psi_j = \frac{\pi 2^{1-4\lambda} \Gamma(j+2\lambda)}{j!(j+\lambda)(\Gamma(\lambda))^2}$

For more properties of Ultraspherical polynomials see [3] by similar way prove (7)

Any function $v(x)$, square integrable in $[0,1]$ may be expressed in terms of Monic Shifted Ultraspherical polynomials as

$$v(x) = \sum_{j=0}^{\infty} a_j \check{\mathbb{C}}_j^{(\lambda)}(x), \quad (8)$$

where the coefficients $a_j$ are given by

$$a_j = \frac{1}{\psi_j} \int_0^1 (x-x^2)^{\lambda-\frac{1}{2}} u(x) \check{\mathbb{C}}_j^{(\lambda)} dx, \quad j = 0,1,\cdots \quad (9)$$

### 3. *The fractional derivatives of $\check{\mathbb{C}}_j^{(\lambda)}(x)$*

This section aims to prove the following theorem for the fractional derivatives of the shifted Monic Ultraspherical polynomials.

**Lemma 3.1** Let $\check{\mathbb{C}}_j^{(\lambda)}(x)$ be a shifted monic Ultraspherical polynomial then
$$D^q \check{\mathbb{C}}_j^{(\lambda)}(x) = 0, \quad j = 0,1,\cdots, \lceil q \rceil - 1, \quad q = 0. \quad (10)$$

**Proof.** This lemma can be easily proved by using relations (2)–(3) with relation (5).

**Theorem 3.2**. The fractional derivative of order $q$ in the Caputo sense for the shifted monic Ultraspherical polynomials is given by

$$D^q \check{\mathbb{C}}_j^{(\lambda)}(x) = \mathcal{H}_j \sum_{r=0}^{[(j-q)/2]} \sum_{k=\lceil q \rceil}^{j-2r} (-1)^{j-r-k} \frac{\Gamma(j-r-\lambda)}{\Gamma(\lambda)\Gamma(r+1)} \frac{2^{j-2r-k}}{\Gamma(k-q+1)\Gamma(j-2r-k+1)} x^{k-q},$$

$$k \geq \lceil q \rceil, \quad j = \lceil q \rceil, \lceil q \rceil + 1, \cdots \quad (11)$$

**Proof.** The analytic form of the MSUPs $\check{\mathbb{C}}_j^{(\lambda)}(x)$ of degree $j$ is given by (5). Using Eqs. (2), (3) and (5) we have

$$D^q \check{\mathbb{C}}_j^{(\lambda)}(x) = \mathcal{H}_j \sum_{r=0}^{[j/2]} \sum_{k=0}^{j-2r} (-1)^{j-r-k} \frac{\Gamma(j-r+\lambda)}{\Gamma(\lambda)\Gamma(r+1)} \frac{2^{j-2r+k}}{\Gamma(k+1)\Gamma(j-2r-k+1)} D^q x^k$$

$$D^q \check{\mathbb{C}}_j^{(\lambda)}(x) = \mathcal{H}_j \sum_{r=0}^{[(j-q)/2]} \sum_{k=0}^{j-2r} (-1)^{j-r-k} \frac{\Gamma(j-r-\lambda)}{\Gamma(\lambda)\Gamma(r+1)} \frac{2^{j-2r-k}}{\Gamma(k-q+1)\Gamma(j-2r-k+1)} x^{k-q},$$

$k \geq \lceil q \rceil$, $j = \lceil q \rceil, \lceil q \rceil + 1, \cdots$

## 4. Method of Solution

This part contains information about the method to solve some types of linear or nonlinear ODEs.

Consider the function $F(x)$ is an integrable square function on $[0,1]$. So, it can be as a form of series in terms of MGPs:

$$F(x) = \sum_{j=0}^{\infty} a_j \check{\mathbb{C}}_N^{(\lambda)}(x), \tag{12}$$

where $a_j$ are arbitrary constant.

Now, assume that the linear DEs as:

$$D^q F(x) + \sum_{i=0}^{r-1} \rho_i(x) D^{s_i} F(x) + g(x) F(x) = G(x). \tag{13}$$

Such as the boundary conditions

$$F^{(k)}(0) = d_k, \quad F^{(l)}(1) = e_l, \quad 0 \leq i, k \leq q - 1. \tag{14}$$

Let $r, s_i, d_k, e_l, k,$ and $i$ are integer number, $s_i < q$, and $i, k, r, s_i \geq 0$, $q$ is fractional number, and $q > 0$. Also, $\rho_i(x), g(x),$ and $G(x)$ are function of $x$.

Hence, we define
$$\begin{aligned} x &= (x_0, x_1, \ldots, x_N)^T, \\ S_N(I) &= span\{\check{\mathbb{C}}_0^{(\lambda)}(x), \check{\mathbb{C}}_1^{(\lambda)}(x), \cdots, \check{\mathbb{C}}_N^{(\lambda)}(x)\}, \\ \rho_i(x) &= (\rho_i(x_0), \rho_i(x_1), \ldots, \rho_i(x_N))^T, \\ g(x) &= (g(x_0), g(x_1), \ldots, g(x_N))^T, \end{aligned} \tag{15}$$

Dente that $a_j$ in equation (8) can be $\boldsymbol{a} = (a_0, a_1, \ldots, a_N)^T$.

Now, by using equations (7) and (8) to get:

$$D^q F(x) = F^{(q)}(x) = \frac{d^q}{dx^q} \sum_{j=0}^{N} a_j \check{\mathbb{C}}_N^{(\lambda)}(x), \tag{16}$$

At $q = 0$ equation (16) will be same equation (12).

By the same substitution at the boundary condition

$$\sum_{j=0}^{N} a_j \left(\check{\mathbb{C}}_N^{(\lambda)}(0)\right)^k = d_k, \quad \sum_{j=0}^{N} a_j \left(\check{\mathbb{C}}_N^{(\lambda)}(1)\right)^l = e_l, \qquad 0 \leq i, k \leq m - 1.$$

By apply the relations (12) and (8) in (13), to get the following result:

$$D^q \sum_{j=0}^{N} a_j \check{\mathbb{C}}_N^{(\lambda)}(x) + \sum_{i=0}^{r-1} \rho_i(x) D^{s_i} \sum_{j=0}^{N} a_j \check{\mathbb{C}}_N^{(\lambda)}(x) + g(x) \sum_{j=0}^{N} a_j \check{\mathbb{C}}_N^{(\lambda)}(x) = G(x). \qquad (17)$$

Equation (13) will be written as:

$$\sum_{j=0}^{N} a_j \left(\check{\mathbb{C}}_N^{(\lambda)}(x)\right)^q + \sum_{i=0}^{r-1} \rho_i(x) \sum_{j=0}^{N} a_j \left(\check{\mathbb{C}}_N^{(\lambda)}(x)\right)^{s_i} + g(x) \sum_{j=0}^{N} a_j \check{\mathbb{C}}_N^{(\lambda)}(x) = G(x). \qquad (14)$$

We can take $a_j$ a common factor from equation (14)

$$a_j \left[\sum_{j=0}^{N} \left(\check{\mathbb{C}}_N^{(\lambda)}(x)\right)^q + \sum_{i=0}^{r-1} \rho_i(x) \sum_{j=0}^{N} \left(\check{\mathbb{C}}_N^{(\lambda)}(x)\right)^{s_i} + g(x) \sum_{j=0}^{N} \check{\mathbb{C}}_N^{(\lambda)}(x)\right] = G(x). \qquad (15)$$

From equation (11) where $x = x_n, n = 0, 1, \cdots, N$ represent by span of $N$ equation (15) will be:

$$a_j \left[\sum_{j=0}^{N} \left(\check{\mathbb{C}}_N^{(\lambda)}(x_n)\right)^q + \sum_{i=0}^{r-1} \rho_i(x_n) \sum_{j=0}^{N} \left(\check{\mathbb{C}}_N^{(\lambda)}(x_n)\right)^{s_i} + g(x_n) \sum_{j=0}^{N} \check{\mathbb{C}}_N^{(\lambda)}(x_n)\right] = G(x_n). \qquad (16)$$

Now, equations (16) represent a system of equation in $a_j$ unknowns. These equations are being linear or nonlinear equations depend on type of equation (9). Where equations (16) are followed equation (9) in being linear or nonlinear equations. By convert equations (16) to matrix $(N + 1) \times (N + 1)$ can be solved too.

## 5. Applications and results

In this section, we solve some examples by the presented method and compare the numerical results with the exact solutions and some earlier works.

**Example 1**. Consider the equation: Bhrawy *et al.* (2014) .

$$D^q u(x) + u(x) = \frac{\Gamma(q_1+1)}{\Gamma(q_1-q+1)} x^{q+q_1} + x^{q_1}, \qquad (17)$$

$0 < q \leq q_1 \leq 1,$

subject to

$u(0) = 0, \ u'(0) = 0 \qquad (18)$

With the exact solution is given by $u(x) = x^{q_1}$.

Tables (5.1), (5.2) observe maximum absolute errors using SMUG method at different values of $\lambda$, and $N$ for example 1. By comparing the numerical results with some earlier works. For that our method appears high accuracy and efficiency.

**Table 5.1:** The observed absolute error using SMUG for Example (1) with various choices of $q_1 = 0.9, q = 0.5$

| N | $\lambda = 1$ | $\lambda = 0.5$ | $\lambda = 0.49$ | $\lambda = -0.49$ | Bhrawy *et al.* (2014) |
|---|---|---|---|---|---|
| 8 | 0 | 1.3696e-16 | 2.8940e-16 | 2.1390e-16 | 1.47e-02 |
| 16 | 1.0430e-14 | 4.1956e-15 | 2.0346e-14 | 1.2416e-15 | 9.45 e-03 |
| 24 | 1.1214e-11 | 4.3510e-13 | 2.4006e-12 | 7.5121e-13 | 4.30 e-03 |
| 32 | 2.0614e-04 | 2.0614e-04 | 2.0614e-04 | 2.0614e-04 | 2.52 e-03 |
| 40 | 6.7773e-08 | 4.9716e-09 | 1.0566e-08 | 7.0142e-09 | 1.44 e-03 |

**Table 5.2:** The observed absolute error using SMUG for Example (1) with various choices of $q_1 = 0.75, q = 0.75$

| N | $\lambda = 1$ | $\lambda = 0.5$ | $\lambda = 0.49$ | $\lambda = -0.49$ | Bhrawy *et al.* (2014) |
|---|---|---|---|---|---|
| 8 | 0 | 3.6746e-16 | 5.1229e-16 | 1.5411e-16 | 2.62 e-02 |
| 16 | 1.0533e-14 | 2.8327e-14 | 5.9736e-15 | 4.3027e-15 | 1.41 e-02 |
| 24 | 6.8690e-12 | 1.4576e-12 | 3.2355e-13 | 1.5140e-12 | 9.08 e-03 |
| 32 | 5.4611e-03 | 5.342e-03 | 5.4611e-03 | 5.4610e-03 | 7.07 e-03 |

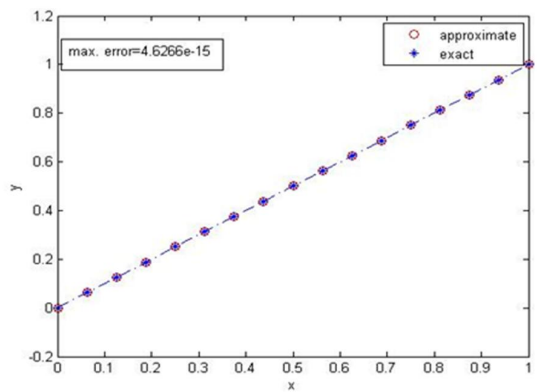

**Fig. 5.1:** The observed absolute error using SMUG for Example (1) with various choices of $q_1 = 0.999, q = 0.999,$
$\lambda = -0.49$, and $N$=16

**Example 2**. Consider the equation Abd-Elhameed and Youssri, (2014)

$$D^q u(x) = 1 - u(x)^2, \quad 0 < q \leq 1, \tag{19}$$

subject to

$$u(0) = 0, \tag{20}$$

With the exact solution is given by $u(x) = \tanh(x)$.

Example 2 Abd-Elhameed and Youssri, (2014) has a maximum absolute errors 3.80e−14 at *N*=14 but our method maximum absolute errors 2.5192e-16

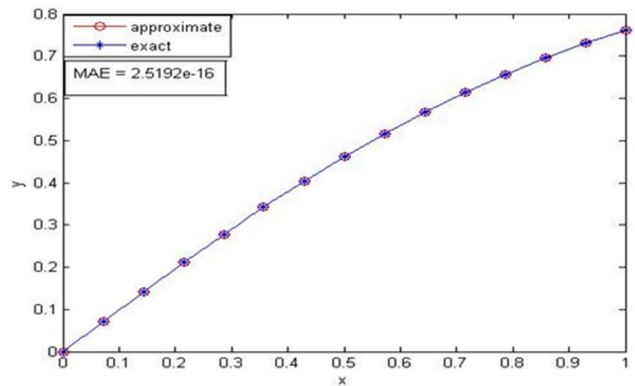

**Fig. 5.2:** The observed absolute error using SMUG for Example (2) with choices of $q = 1$, *N*=14, and $\lambda = -0.49$

**Example 3**. Non-linear fourth order BVP: ( Abdelhakem *et al.*, 2019)

$$16u^{(4)}(x) - 6e^{-4u(x)} = -12(1.5 + 0.5x)^{-4}; \quad -1 \leq x \leq 1, \tag{21}$$

With boundary conditions:
$$\left.\begin{array}{l} u(-1) = 0, \quad u(1) = \ln(2), \\ u'(-1) = 0.5, \quad u'(1) = 0.25, \end{array}\right\} \tag{22}$$
with analytical solution:
$$u(x) = \ln(1.5 + 0.5x). \tag{23}$$

**Table 5.3:** The observed absolute error using SMUG for Example (3)

| N | $\alpha =1$ | $\alpha =0.5$ | $\alpha =-0.49$ | $\alpha =0.49$ | Abdelhakem *et al.*, (2019) |
|---|---|---|---|---|---|
| 8 | 6.9315e-06 | 6.9315e-06 | 6.9315e-06 | 6.9315e-06 | - |
| 10 | 6.5393e-08 | 6.5393e-08 | 6.5393e-08 | 6.5393e-08 | 7.5991e-05 |
| 12 | 1.9845e-09 | 1.9845e-09 | 1.9845e-09 | 1.9845e-09 | 3.5789e-07 |
| 14 | 2.2012e-11 | 2.2018e-11 | 2.2006e-11 | 2.2010e-11 | 1.2290e-09 |
| 16 | 1.2683e-13 | 1.3029e-13 | 1.2684e-13 | 1.2785e-13 | 3.5426e-12 |
| 18 | 6.4429e-14 | 6.4399e-14 | 6.4114e-14 | 6.4494e-14 | 1.5099e-14 |
| 20 | 4.2415e-16 | 4.3027e-16 | 4.2409e-16 | 4.2418e-16 | 4.4409e-16 |

**Conclusion**

An efficient and accurate numerical algorithm based on the shifted Monic Ultraspherical Galerkin spectral method is proposed for solving the FDEs. The problem is reduced to the solution of system of simultaneous nonlinear algebraic equations. To the best of our knowledge, this is the first work concerning the SMU Galerkin spectral method algorithm for solving general fraction order differential equations. Numerical examples were given to demonstrate the validity and applicability of the algorithm. The results show that the method is simple and accurate. In fact, by selecting few collocation points, excellent numerical results are obtained.